\documentclass[12pt]{article}
\usepackage{amssymb,amsmath, amsthm, amscd,ifthen}
\usepackage[dvips]{graphicx}
\usepackage{indentfirst}

\renewcommand{\le}{\leqslant}
\renewcommand{\ge}{\geqslant}

\begin{document}

\title{On hypergraph cliques with chromatic number 3 and a given number of vertices\footnote{This work was supported by 
the grant of RFBR N 12-01-00683, the grant of the Russian President N MD-8390.2010.1, and the grant NSh-8784.2010.1.}}

\author{D.D. Cherkashin, A.B. Kulikov, A.M. Raigorodskii}

\date{}

\maketitle

\begin{abstract}

In 1973 P. Erd\H{o}s and L. Lov\'asz noticed that any hypergraph whose edges are pairwise intersecting has chromatic number 2 or 3. In 
the first case, such hypergraph may have any number of edges. However, Erd\H{o}s and Lov\'asz proved that in the second case, the number 
of edges is bounded from above. For example, if a hypergraph is $ n $-uniform, has pairwise intersecting edges, and has chromatic number 
3, then the number of its edges does not exceed $ n^n $. Recently D.D. Cherkashin improved this bound (see \cite{Ch}). In this paper, we further 
improve it in the case when the number of vertices of an $n$-uniform hypergraph is bounded from above by $ n^m $ with some $ m = m(n) $. 

\end{abstract}

\section{Introduction and formulation of the main result}

This work is devoted to a problem in extremal hypergraph theory, which goes back to P. Erd\H{o}s and L. Lov\'asz (see \cite{EL}). Before 
giving an exact statement of the problem, we recall some definitions and introduce some notation. 

Let $ H = (V,E) $ be a hypergraph without multiple edges. We call it {\it $ n $-uniform}, if any of its edges has cardinality $ n $: for every
$ e \in E $, we have $ |e| = n $. By the {\it chromatic number} of a hypergraph $ H = (V,E) $ we mean the minimum number 
$ \chi(H) $ of colors needed to color 
all the vertices in $ V $ so that any edge $ e \in E $ contains at least two vertices of some different colors. Finally, a hypergraph is said 
to form a {\it clique}, if its edges are pairwise intersecting. 

In 1973 Erd\H{o}s and Lov\'asz noticed that if an $ n $-uniform hypergraph $ H = (V,E) $ forms a clique, then $ \chi(H) \in \{2,3\} $. They also 
observed that in the case of $ \chi(H) = 3 $, one certainly has $ |E| \le n^n $ (see \cite{EL}). Thus, the following definition has been 
motivated: 
$$
M(n) = \max \{|E|: ~ \exists ~ {\rm an} ~ n-{\rm uniform} ~ {\rm clique} ~ H = (V,E) ~ {\rm with} ~ \chi(H) = 3\}.
$$
Obviously such definition has no sense in the case of $ \chi(H) = 2 $. 

\vskip+0.2cm

\noindent{\bf Theorem 1 (P. Erd\H{o}s, L. Lov\'asz, \cite{EL}).} {\it The inequalities hold 
$$
n! \left(\frac{1}{1!} + \frac{1}{2!} + \ldots + \frac{1}{n!}\right) \le M(n) \le n^n.
$$}

\vskip+0.2cm

Almost nothing better has been done during the last 35 years. In the book \cite{TJ} the estimate $ M(n) \le \left(1-\frac{1}{e}\right) n^n $ 
is mentioned as ``to appear''. However, we have not succeeded in finding the corresponding paper. 

At the same time, another quantity $ r(n) $ was introduced in \cite{Lov}:
$$
r(n) = \max \{|E|: ~ \exists ~ {\rm an} ~ n-{\rm uniform} ~ {\rm clique} ~ H = (V,E) ~ {\rm s.t.} ~ \tau(H) = n\},
$$
where $ \tau(H) $ is the {\it covering number} of $ H $, i.e., 
$$
\tau(H) = \min \{|f|: ~ f \subset V, ~ \forall ~ e \in E ~~ f \cap e \neq \emptyset\}.
$$
Clearly, for any $ n $-uniform clique $ H $, we have $ \tau(H) \le n $ (since every edge forms a cover), and if $ \chi(H) = 3 $, then 
$ \tau(H) = n $. Thus, $ M(n) \le r(n) $. Lov\'asz noticed that for $ r(n) $ the same estimates as in Theorem 1 apply and conjectured that 
the lower estimate is best possible. In 1996 P. Frankl, K. Ota, and N. Tokushige (see \cite{FOT}) disproved this conjecture and showed that 
$ r(n) \ge \left(\frac{n}{2}\right)^{n-1} $. 

In \cite{Ch} D.D. Cherkashin discovered a new upper bound for the initial value $ M(n) $ which is actually true for $ r(n) $ as well. 

\vskip+0.2cm

\noindent{\bf Theorem 2 (D.D. Cherkashin, \cite{Ch}).} {\it There exists a constant $ c > 0 $ such that
$$
M(n) \le c n^{n-\frac{1}{2}} \ln n.
$$}

\vskip+0.2cm

To formulate the main result of this paper we take any natural numbers $ n $, $ m \ge 2 $ and put $ q(n,m) = \left[\frac{n}{2m}\right] $, 
$$ 
A(n,m) = \sum_{i=0}^{2q(n,m)} {n^m \choose i}. 
$$
We note that of course 
$$
A(n,m) \le \left(\frac{n}{m}+1\right) {n^m \choose 2q(n,m)} \le \left(\frac{n}{m}+1\right) 
\left(\frac{e n^m}{2q(n,m)}\right)^{n/m} = n^n \cdot A'(n,m), 
$$
where
$$
A'(n,m) = \left(\frac{n}{m}+1\right) \left(\frac{e}{2q(n,m)}\right)^{n/m}.
$$
Obviously, if $ m $ is a function of $ n $, which is $ o(n) $ as $ n \to \infty $, then 
$$ 
A'(n,m) = \frac{m}{n \omega(n)},
$$
where $ \omega(n) \to \infty $ as $ n \to \infty $. Thus, $ A(n,m) = o\left(m n^{n-1}\right) $. 

\vskip+0.2cm

\noindent{\bf Theorem 3.} {\it Let $ m \ge 2 $ be any function of $ n \in {\mathbb N} $ which is $ o(n) $ as $ n \to \infty $; 
moreover, $ m(n) \le \frac{n}{2} $. For any 
$ n \ge 4 $ and any $ n $-uniform clique $ H = (V,E) $ with $ \chi(H) = 3 $ and $ |V| \le n^{m(n)} $, we have  
$$
|E| \le 4m(n) n^{n-1} + A(n,m(n)) = (4+o(1))m(n) n^{n-1}.
$$}

\vskip+0.2cm

Clearly, if $ m(n) \le c \sqrt{n}\ln n $ with some constant $ c > 0 $, then the bound in Theorem 3 is stronger than the bound in 
Theorem 2. Note that the number of vertices in any $ n $-uniform clique with chromatic number 3 does not exceed $ 4^n $ (see \cite{EL}).
Unfortunately, $ n^{\sqrt{n}\ln n} = e^{o(n)} $, so that Theorem 3 does not cover all possible values of $ |V| $.

\section{Proof of Theorem 3}

Fix an $ n \ge 4 $ and put $ m = m(n) $, $ q = q(n,m) $, $ A = A(n,m) $. Fix an $ n $-uniform clique $ H = (V,E) $ with 
$ \chi(H) = 3 $ and $ |V| \le n^m $. For any set $ W \subseteq V $, denote by $ E(W) $ the set of all edges $ B \in E $ such that 
$ W \subseteq B $. Also denote by $ E_W $ the set of all edges $ B \in E $ such that $ W \cap B \neq \emptyset $. Clearly $ E(W) \subseteq 
E_W $. Let 
$$
Q = \{1,2,3, \dots, q\} \cup \{n-q+1, n-q+2, \dots, n\}.
$$
The two parts which form the set $ Q $ do not intersect and do not cover the whole set $ \{1, \dots, n\} $, since $ m \ge 2 $. Moreover, 
$ Q $ is not empty, since $ m \le \frac{n}{2} $ and so $ q \ge 1 $. 

\vskip+0.2cm

\noindent{\bf Lemma 1.} {\it Let $ W \subseteq V $, $ i = |W| $. Either there exists a vertex $ x \in W $ such that $ {\rm deg}~x \ge 
\frac{|E|-A}{i} $, or there exist two edges $ B_1, B_2 \in E $ such that $ B_1, B_2 \not \in E_W $ and $ |B_1 \cap B_2| \not \in Q $.}

\vskip+0.2cm

\paragraph{Proof of Lemma 1.} If there exists a vertex $ x \in W $ such that $ {\rm deg}~x \ge \frac{|E|-A}{i} $, then we are done. 
If there are no such vertices, then 
$$
|E_W| \le \sum_{x \in W} {\rm deg}~x < |E| - A.
$$
Therefore, $ |E \setminus E_W| > A $. We have to show that there exist two edges $ B_1, B_2 \in E \setminus E_W $ with $ |B_1 \cap B_2| \not \in Q $.
Suppose to the contrary that for any $ B_1, B_2 \in E \setminus E_W $, we have $ |B_1 \cap B_2| \in Q $. We shall consequently prove that 
$ |E \setminus E_W| \le A $ obtaining a contradiction and thus completing the proof of Lemma 1.

In principle, it is possible just to cite the paper \cite{F}. 
We use instead a version of the linear algebra method in combinatorics (see \cite{BF} and \cite{Rai1}). To any edge $ B $ from 
$ E \setminus E_W $ we assign a vector $ {\bf x} = (x_1, \dots, x_v) \in \{0,1\}^v $, where $ v = |V| \le n^m $ and $ x_{\nu} = 1 $, if and 
only if $ \nu \in B $. In particular, $ x_1 + \ldots + x_v = n $. Let $ E \setminus E_W \to \{{\bf x}_1, \dots, {\bf x}_s\} $. 

Denote by $ ({\bf x},{\bf y}) $ the Euclidean inner product of vectors $ {\bf x}, {\bf y} $. Note that if $ B, B' \in E \setminus E_W $ and 
$ {\bf x}, {\bf x}' $ are the corresponding vectors, then $ |B \cap B'| = ({\bf x},{\bf x}') $.

Take an arbitrary vector $ {\bf x}_{\nu} $, $ \nu \in \{1, \dots, s\} $, and consider the polynomial
$$
F_{{\bf x}_{\nu}}({\bf y}) = \prod_{j \in Q \setminus \{n\}} (j - ({\bf x}_{\nu}, {\bf y})) \in {\mathbb R}[y_1, \dots, y_v].
$$
Eventually, we get $ s $ polynomials $ F_{{\bf x}_{1}}, \dots, F_{{\bf x}_{s}} $. All of them depend on $ v $ variables and have degree not 
exceeding the quantity $ |Q| \le 2q $. Of course any such polynomial is a linear combination of some monomials which are of type 
$$
1, ~~~ y_{\nu_1}^{\alpha_{\nu_1}} \cdot \ldots \cdot y_{\nu_r}^{\alpha_{\nu_r}}, ~~~ \alpha_{\nu_1}, \dots, \alpha_{\nu_r} \ge 1, ~~~ 
\alpha_{\nu_1} + \ldots + \alpha_{\nu_r} \le |Q| \le 2q.
$$
Replace each monomial of this type by $ y_{\nu_1} \cdot \ldots \cdot y_{\nu_r} $. Denote by 
$ F'_{{\bf x}_{1}}, \dots, F'_{{\bf x}_{s}} $ the resulting polynomials. They also depend on $ v $ variables and have degree not 
exceeding the quantity $ |Q| \le 2q $. Moreover, they span a linear space whose dimension is less then or equal to  
$$
\sum_{r=0}^{2q} {v \choose r} \le \sum_{r=0}^{2q} {n^m \choose r} = A.
$$
At the same time $ F'_{{\bf x}_{\nu}}({\bf y}) = F_{{\bf x}_{\nu}}({\bf y}) $, provided $ {\bf y} \in \{0,1\}^v $ and 
$ \nu \in \{1, \dots, s\} $. 

To show that $ s = |E \setminus E_W| \le A $ (which we need to complete the proof) it suffices to establish the linear independence of the polynomials 
$ F'_{{\bf x}_{1}}, \dots, F'_{{\bf x}_{s}} $ over $ {\mathbb R} $. Assume that 
$$
c_1 F'_{{\bf x}_{1}}({\bf y}) + \ldots + c_s F'_{{\bf x}_{s}}({\bf y}) = 0.
$$
Let $ {\bf y} = {\bf x}_{\nu} $, $ \nu \in \{1, \dots, s\} $. Then $ ({\bf x}_{\nu},{\bf y}) = ({\bf x}_{\nu},{\bf x}_{\nu}) = n $ and 
$$ 
F'_{{\bf x}_{\nu}}({\bf y}) = F_{{\bf x}_{\nu}}({\bf y}) = F_{{\bf x}_{\nu}}({\bf x}_{\nu}) \neq 0. 
$$
However, if $ \mu \neq \nu $, then $ ({\bf x}_{\mu},{\bf y}) = ({\bf x}_{\mu},{\bf x}_{\nu}) \in Q \setminus \{n\} $, that is, 
$$ 
F'_{{\bf x}_{\mu}}({\bf y}) = F_{{\bf x}_{\mu}}({\bf y}) = F_{{\bf x}_{\mu}}({\bf x}_{\nu}) = 0. 
$$
Therefore, $ c_{\nu} = 0 $ for every $ \nu $, and we are done. Lemma 1 is proved. 
 
\vskip+0.2cm

\noindent{\bf Lemma 2.} {\it Let $ W \subseteq V $, $ i = |W| $, $ j = |E(W)| $. Assume that there exist two edges 
$ B_1, B_2 \in E \setminus E_W $ such that $ |B_1 \cap B_2| \not \in Q $. Put $ \tau = 1 + \frac{1}{4m} $. 
Either there exists an $ x \not \in W $ such that $ |E(W \cup \{x\})| \ge \frac{j\tau}{n} $, or there exist $ x,y \not \in W $ such that 
$ |E(W \cup \{x,y\})| \ge \frac{j\tau^2}{n^2} $.}


\paragraph{Proof of Lemma 2.} Let $ l = |B_1 \cap B_2| \not \in Q $. Consider the set $ E(W) $. Since $ H $ is a clique, any edge 
$ B \in E(W) $ intersects both $ B_1 $ and $ B_2 $. Either $ B $ intersects the set $ B_1 \cap B_2 $, or it has common vertices with 
both $ B_1 \setminus (B_1 \cap B_2) $ and $ B_2 \setminus (B_1 \cap B_2) $. Denote by $ E_1 $ the set of edges of the first type; 
$ E_2 = E(W) \setminus E_1 $. By pigeon-hole principle, there is an $ x \in B_1 \cap B_2 $ such that $ x $ belongs to at least 
$ \frac{|E_1|}{l} $ edges from $ E_1 $; also there are $ x \in B_1 \setminus (B_1 \cap B_2) $ and $ y \in B_2 \setminus (B_1 \cap B_2) $ such 
that the set $ \{x,y\} $ belongs to at least $ \frac{|E_2|}{(n-l)^2} $ edges from $ E_2 $. It remains to show that for any partition 
$ E(W) = E_1 \cup E_2 $, we have 
$$
{\rm either} ~~~ \frac{|E_1|}{l} \ge \frac{j\tau}{n}, ~~~ {\rm or} ~~~ \frac{|E_2|}{(n-l)^2} \ge \frac{j\tau^2}{n^2},
$$
which is equivalent to
$$
\max \left\{\frac{|E_1|^2}{j^2l^2}, \frac{|E_2|}{j(n-l)^2}\right\} \ge \frac{\tau^2}{n^2}.
$$

Here the worst case is that of $ \frac{|E_1|^2}{j^2l^2} = \frac{|E_2|}{j(n-l)^2} $. Let $ a = |E_1| $. Then $ |E_2| = j-a $ and we have
$ \frac{a^2}{j^2l^2} = \frac{j-a}{j(n-l)^2} $. Solving this equation we get
$$
a = \frac{jl^2 + \sqrt{(jl^2)^2+4j^2l^2(n-l)^2}}{2(n-l)^2}.
$$
Of course the value of $ |E_1| $ (which is integer) may differ from the real number $ a $. However, we do certainly know that 
$$
\max \left\{\frac{|E_1|^2}{j^2l^2}, \frac{|E_2|}{j(n-l)^2}\right\} \ge \frac{a^2}{j^2l^2}.
$$
Thus, we need to prove that $ \frac{a}{jl} \ge \frac{\tau}{n} $ or that $ \frac{an}{jl} \ge \tau $. We have
$$
\frac{an}{jl} = \frac{l + \sqrt{l^2 + 4(n-l)^2}}{2(n-l)^2} n = \frac{ln}{2(n-l)^2} + 
\sqrt{\left(\frac{ln}{2(n-l)^2}\right)^2 + \frac{n^2}{(n-l)^2}} \ge 1 + \frac{ln}{2(n-l)^2}.
$$

The function $ \frac{ln}{2(n-l)^2} $ is monotone increasing in $ l $. Since $ l \not \in Q $, we may use the bound 
$ l \ge \frac{n}{2m} $. Consequently, 
$$
\frac{an}{jl} \ge 1 + \frac{ln}{2(n-l)^2} \ge 1+ \frac{n^2}{4m (n-l)^2} \ge 1 + \frac{1}{4m} = \tau.
$$

Lemma 2 is proved. 

\vskip+0.2cm

\paragraph{Completion of the proof of Theorem 3.}

Let 
$$
k = \min \left\{|W|: ~ W \subseteq V, ~ \exists ~ x \in W ~~ {\rm deg}~x \ge \frac{|E|-A}{|W|}\right\}.
$$
The quantity $ k $ is well-defined. Indeed, take any edge $ W \in E $. Since $ H $ is a clique, $ W $ intersects all the edges from $ E $ and 
so there exists an $ x \in W $ with $ {\rm deg}~x \ge \frac{|E|}{|W|} > \frac{|E|-A}{|W|} $. 

Let $ W_0 $ be a set on which the value of $ k $ is attained. Take a vertex $ x \in W_0 $ that has $ {\rm deg}~x \ge \frac{|E|-A}{k} $. The last 
inequality can be rewritten as $ |E(\{x\})| \ge \frac{|E|-A}{k} $. If $ k \ge 2 $, we may apply Lemmas 1 and 2 to $ W = \{x\} $. Thus, we obtain 
either a set $ W' $ of two elements with $ |E(W')| \ge \frac{|E|-A}{k} \cdot \frac{\tau}{n} $ or a set $ W'' $ of three elements with 
$ |E(W'')| \ge \frac{|E|-A}{k} \cdot \frac{\tau^2}{n^2} $. We continue this process until we get a set $ W $ with $ |W| = k $ and 
$ |E(W)| \ge \frac{|E|-A}{k} \cdot \frac{\tau^{k-1}}{n^{k-1}} $ (even if $ k = 1 $, we do have such a set). 

In \cite{EL} Erd\H{o}s and Lov\'asz proved that for any $ n $-uniform clique $ H = (V,E) $ with chromatic number 3, if $ W \subseteq V $ is of 
cardinality $ k $, then $ |E(W)| \le n^{n-k} $. In our case, we have $ \frac{|E|-A}{k} \cdot \frac{\tau^{k-1}}{n^{k-1}} \le n^{n-k} $. 
Therefore, 
$$
|E| \le k \cdot n^{n-k} \cdot \frac{n^{k-1}}{\tau^{k-1}} + A = k \frac{n^{n-1}}{\tau^{k-1}} + A. 
$$

To complete the proof of Theorem 3 it remains to show that for any $ k $, $ \frac{k}{\tau^{k-1}} \le 4m $. It is very easy to see that the 
maximum value of the quantity $ \frac{k}{\tau^{k-1}} $ is attained on $ k = 4m $, and we are done. 

\section{A refinement in Theorem 3}

For $ m = 2 $, one can prove a simple result, which is however substantially better than that of Theorem 3. 

\vskip+0.2cm

\noindent{\bf Theorem 4.} {\it Let $ H = (V,E) $ be any $ n $-uniform clique with $ \chi(H) = 3 $. Put $ v = |V| $. Assume that $ v \le 
\frac{n^2}{c} $, where $ c $ may be any function of $ n $ such that $ c(n) \in (1,n) $. Put 
$$ 
d = c e^{\frac{2}{ec}-1}.
$$
Then 
$$ 
|E| \le (1+o(1)) \frac{e^{3/2}}{\sqrt{c}} (n/d)^n. 
$$} 

\vskip+0.2cm

If $ c $ is a constant, then we get an exponential improvement for the Erd\H{o}s and Lov\'asz bound by $ n^n $. Otherwise, the improvement is 
even more considerable. 

\paragraph{Proof of Theorem 4.} Take an arbitrary integer $ a \in (1,n) $ and consider all the $ a $-element subsets of $ V $. The number of 
such subsets is $ {v \choose a} $. On the one hand, any edge from $ E $ contains exactly $ {n \choose a} $ subsets. On the other hand, any 
subset is contained in at most $ n^{n-a} $ edges (see \cite{EL}). So the number of edges does not exceed the quantity 
$ \frac{n^{n-a} {v \choose a}}{{n \choose a}} $. To estimate this quantity we use the bound $ {v \choose a} \le \frac{v^a}{a!} $ and 
the Stirling formula. Hence, 
$$
\frac{n^{n-a} {v \choose a}}{{n \choose a}} \le \frac{n^{n-a} v^a}{a! {n \choose a}} \le \frac{n^{n+a}}{c^a \frac{n!}{(n-a)!}}.
$$

Now put 
$$ 
a = \left[\left(1-\frac{1}{ec}\right)n\right] + 1. 
$$
Then $ n - a \le \frac{n}{ec} $, so that 
$$
\frac{n!}{(n-a)!} \sim \frac{\sqrt{2\pi n} \left(\frac{n}{e}\right)^n}{\sqrt{2\pi (n-a)} \left(\frac{n-a}{e}\right)^{n-a}} \ge 
\sqrt{ec} \cdot n^a (e^2c)^{n-a} e^{-n}
$$
and
$$
|E| \le (1+o(1)) \frac{n^{n+a}}{c^a \sqrt{ec} \cdot n^a (e^2c)^{n-a} e^{-n}} = (1+o(1)) \frac{n^n}{\sqrt{ec} \cdot c^n e^{-n} e^{2n-2a}} \le
$$
$$ 
\le (1+o(1)) \frac{n^n}{\sqrt{ec} \cdot c^n e^{-n} e^{\frac{2n}{ec}-2}} = (1+o(1)) \frac{e^{3/2}}{\sqrt{c}} (n/d)^n. 
$$

Theorem 4 is proved. 

\vskip+0.2cm

Note that for constant values of $ c $, the choice of $ a $ in the proof was nearly optimal.

\end{document}